\documentclass{amsart}

\theoremstyle{plain}
\newtheorem{Thm}{Theorem}[section]
\newtheorem{Lem}[Thm]{Lemma}
\newtheorem{Prop}[Thm]{Proposition}
\newtheorem{Cor}[Thm]{Corollary}

\theoremstyle{definition}
\newtheorem{Def}[Thm]{Definition}

\theoremstyle{remark}
\newtheorem{remark}[Thm]{Remark}
\numberwithin{equation}{section}

\def\CC{{\mathbb C}}
\def\GH{{\mathcal H}}
\def\LH{{\mathcal L}(\GH)}

\def\DD{{\mathbb D}}
\def\TT{{\mathbb T}}

\def\f1{\mathbb{1}}
\def\HH{{\mathcal H}}
\def\LL{{\mathcal L}}
\def\EE{{\mathcal E}}
\def\DT{{\mathcal D_T}}
\def\DTT{{\mathcal D_{T^*}}}
\newcommand{\TTT}{\mathbf T}
\newcommand{\KK}{\mathfrak K}
\newcommand{\Dcal}{\mathcal D}

\def\dfrac{\displaystyle\frac}

\def\Dim{\mathop{\rm dim}\nolimits}

\begin{document}

\title[Complex symmetric contractions]{The characteristic function of a complex symmetric contraction}
\author{Nicolas Chevrot} 

\address{Institut Camille Jordan, UFR de Math\'ematiques,
Universit\'e Claude Bernard Lyon I, 69622 Villeurbanne Cedex,
France.}
\email{chevrot@math.univ-lyon1.fr}

\author{Emmanuel Fricain}

\address{Institut Camille Jordan, UFR de Math\'ematiques,
Universit\'e Claude Bernard Lyon I, 69622 Villeurbanne Cedex,
France.}
\email{fricain@math.univ-lyon1.fr}

\author{Dan Timotin}

\address{Institute of Mathematics of the Romanian
Academy, P.O. Box 1-764, Bucharest 014700, Romania.}
\email{Dan.Timotin@imar.ro}

\thanks{}



\begin{abstract}
It is shown that a contraction on a Hilbert space is complex symmetric if and
only if the values of its characteristic function are all symmetric with 
respect to a fixed conjugation. Applications are given to the description
of complex symmetric contractions with defect indices equal to~2.
\end{abstract}
\bibliographystyle{amsplain}

\maketitle

\section{Introduction}

Complex symmetric operators on a complex Hilbert space are characterized by
the existence of an orthonormal basis with respect to which their
matrix is symmetric. Their theory is therefore connected
with the theory of symmetric matrices, which is a classical topic in linear
algebra. A more intrinsic definition implies the introduction of a conjugation
in the Hilbert space, that is, an antilinear, isometric and involutive map,
with respect to which the symmetry is defined.
Such operators or matrices apppear naturally in many different areas of mathematics
and physics; we refer to~\cite{GP} for more about the history of the subject
and its connections to other domains, as well as for an extended list of references.
 
The interest in complex symmetric operators has been recently revived by
the work of Garcia and Putinar~\cite{Garcia1, Garcia2, GP}. In their papers
a general framework is established for such operators, and it is shown
that large classes of operators on a Hilbert space can be studied 
in this framework. The examples are rather diverse: normal operators
are complex symmetric, for instance, but also certain types of Volterra
and Toeplitz
operators, as well as the so-called compressed shift on the
functional model spaces $H^2\ominus \phi H^2$, where $\phi$ denotes 
a nonconstant inner function.

The purpose of this paper is to explore further the generalizations 
of this last example. The natural context is the model theory of completely
non unitary 
contractions developed by Sz. Nagy and Foias~\cite{SNF}. The main result
is a criterium for a contraction to be complex symmetric in terms
of its characteristic function. In the sequel some applications
of this result are given.

The plan of the paper is the following. The next section presents preliminary
material. Section~3 contains the announced criterium. In Section 4 one discusses
$2\times 2$ inner characteristic functions, and the results are applied in the last section 
in order to obtain a series of examples of complex symmetric contractions with
defect indices~2.

\section{Preliminaries}

\subsection{Complex symmetric operators}
We first recall some basic  facts from~\cite{Garcia1, Garcia2, GP}. 
Let $\HH$ be a complex Hilbert space, and $\LL(\HH)$ the algebra of all bounded linear 
operators on $\HH$. A \emph{conjugation} $C$ on $\HH$ is an anti-linear, isometric and involutive map; thus $C^2=I$, and $\langle Cf,Cg \rangle=\langle g,f\rangle$ for all $f,g\in \HH$. It is easy to see that if $C$ is a conjugation, and $S\in\LL(\HH)$ is a symmetry (that is, a unitary involutive
operator), then $CS$ is also a conjugation. 


For a fixed conjugation operator $C$ on $\HH$, we say that a linear operator $T$ on $\HH$ is \emph{$C$-symmetric} if $T=CT^*C$. One sees
immediately that if $T$ is $C$-symmetric, then $T^*$ is $C$-symmetric.
Then $T\in\LL(\HH)$ is called \emph{complex symmetric} if there exists a conjugation $C$ on $\HH$ such that $T$ is $C$-symmetric. 
Among  various examples of complex symmetric operators \cite{GP},
we mention the class of normal operators; in particular, unitary 
operators are complex symmetric. Also, direct sums of complex symmetric operators
are complex symmetric.

Complex symmetric operators can also be characterized in terms of certain matrix representations,
as shown by the following result from~\cite{GP}.

\begin{Lem}\label{Prop:symmetrie}
Let  $C$ be a conjugation on $\HH$. Then:
\begin{enumerate}
\item[$\mathrm{(i)}$] There exists an orthonormal basis $(e_n)_{n=1}^{\Dim\HH}$ of $\HH$ such that $Ce_n=e_n$ for all $n$; such a basis is called a \emph{$C$-real orthonormal basis} for $\HH$.
\item[$\mathrm{(ii)}$]  $T\in\LL(\HH)$ is $C$-symmetric if an only if there exists a $C$-real orthonormal basis $(e_n)_{n=1}^{\Dim\HH}$ for $\HH$ such that 
\[
\langle Te_n,e_m\rangle=\langle Te_m,e_n\rangle,\qquad\forall n,m\geq 1.
\]
\end{enumerate}
\end{Lem}

\subsection{Characteristic functions and model operators}
The characteristic function for a contraction and the construction of the basic functional model is developed by B. Sz.-Nagy and C. Foias~\cite{SNF}, which is the main source 
for this subsection.  Let $T\in\LH$ be a contraction, that is, $\|T\|\le 1$. 
There is a unique decomposition $\HH=\HH_0\oplus \HH_u$
such that $T\HH_0\subset \HH_0$, $T\HH_u\subset \HH_u$ and $T_{|\HH_u}$ is unitary, whereas $T_{|\HH_0}$ is \emph{completely non-unitary} (c.n.u.), that is, $T_{|\HH_0}$ is not unitary on any of its invariant subspaces.

The operator $D_T=(I-T^*T)^{1/2}$ is called the \emph{defect operator} of $T$. The \emph{defect spaces} of $T$ are $\mathcal D_T=\overline{D_T\HH}$, $\mathcal D_{T^*}=\overline{D_{T^*}\HH}$,
and the \emph{defect indices} $\partial_T=\dim\DT$, $\partial_{T^*}=\dim\DTT$.
Since $D_T=D_{T_0}\oplus 0$, $D_{T^*}=D_{T_0^*}\oplus 0$, we have $\mathcal D_T=\mathcal D_{T_0}$ and $\mathcal D_{T^*}=\mathcal D_{T_0^*}$. 

We say that $T\in C_{0.}$ if $T^n\to 0$ strongly, and $T\in C_{.0}$ if $T^*\in C_{0.}$; 
also, $C_{00}=C_{0.}\cap C_{.0}$.

Suppose $\EE, \EE'$ are Hilbert spaces, and  $\Theta:\DD\to\LL(\EE, \EE')$ 
is a contraction-valued analytic function. One can decompose 
$\EE=\EE_p\oplus\EE_u$, 
$\EE'=\EE'_p\oplus\EE'_u$, such that:
\begin{enumerate}
\item[---]
 for all $z\in\DD$, $\Theta(z)\EE_p\subset \EE'_p$,
 $\Theta(z)\EE_u\subset \EE'_u$;
 \item[---] if $\Theta=\Theta_p\oplus \Theta_u$ is the corresponding decomposition
 of~$\Theta$, then $\Theta_p$ is \emph{pure},
that is, $\|\Theta_p(0)h\|<\|h\|$ for all $h\in\EE_p$, $h\not=0$, 
while $\Theta_u$ is a unitary constant.
\end{enumerate}
$\Theta_p$ is then called the pure part of $\Theta$.

One says~\cite{SNF} that two contractive analytic
functions $\Theta:\DD\to\LL(\EE, \EE_*)$, $\Theta':\DD\to\LL(\EE', \EE'_*)$
\emph{coincide} if there are unitaries $U:\EE\to\EE'$, $U_*:\EE_*\to\EE'_*$,
such that $\Theta(z)=U_*^*\Theta'(z) U$ for all $z\in\DD$.

The \emph{characteristic function} of $T$ is an operator valued function $\Theta_T(\lambda):\mathcal D_T\to \mathcal D_{T^*}$ defined for $\lambda\in\DD$ by
\begin{equation}\label{eq:fct-carac1}
\Theta_T(\lambda):=-T+\lambda D_{T^*}(I-\lambda T^*)^{-1}D_T|\mathcal D_T.
\end{equation}
$\Theta_T$ is a pure contraction-valued analytic function on $\DD$, and one sees easily that
$\Theta_T=\Theta_{T_0}$. 

For $\EE$ a Hilbert space, we denote by $L^2(\EE)$ the Lebesgue space of measurable functions
$f:\TT\to\EE$ of square integrable norm, and by $H^2(\EE)\subset L^2(\EE)$
the Hardy space of functions whose negative Fourier coefficients vanish. $P_+$
is the orthogonal projection onto $H^2(\EE)$, and $P_-=I-P_+$.

If we are given an arbitrary  contraction-valued analytic function
$\Theta:\DD\to\LL(\EE,\EE')$ ($\EE, \EE_*$ Hilbert spaces), one defines the
\emph{model space} associated to $\Theta$ by
\begin{equation}\label{eq:espace-modele1}
\mathfrak K_{\Theta}=\left(H^2(\EE_*)\oplus \overline{(I-\Theta^*\Theta)^{1/2} L^2(\EE)}\right)\ominus \{\Theta_T f\oplus (I-\Theta^*\Theta)^{1/2} f:f\in H^2(\EE)\},
\end{equation}
and the \emph{model operator} $\mathbf T_{\Theta}\in\LL(\mathfrak K_{\Theta})$ by 
\begin{equation}\label{eq:operateur-modele}
\mathbf T_{\Theta}(f\oplus g)=P_{\mathfrak K_{\Theta}}(zf\oplus zg)
\end{equation}
($P_{\mathfrak K_{\Theta}}$ is the orthogonal projection onto $\mathfrak K_{\Theta}$).
Then $\mathbf T_{\Theta}$ is a c.n.u. contraction, and its characteristic
function coincides with the pure part of $\Theta$.

If we start with a contraction $T$, and apply the previous constructions
to $\Theta_T$, the resulting operator $\mathbf T_{\Theta_T}$ is unitarily equivalent to $T_0$
(the completely non-unitary part of $T$).


A contractive analytic function $\Theta$ is called \emph{inner} if its boundary
values $\Theta(e^{it})$ are isometries a.e. on $\TTT$. If $T$ is c.n.u., then $T\in C_{.0}$
if and only if $\Theta_T$ is inner.

\section{The main theorem}



Our main result gives a criterion for complex symmetric contractions.

\begin{Thm}\label{Thm:main}
Let $T$ be a contraction on the Hilbert space $\HH$. Then the following are equivalent:
\begin{enumerate}
\item[\rm (i)] $T$ is complex symmetric.
\item[\rm (ii)] There exists an anti-linear map $J:\mathcal D_T\to\mathcal D_{T^*}$ which is isometric, onto and satisfies
\begin{equation}\label{eq:caract-symmetric}
\Theta_T(z)=J\Theta_T(z)^*J,\qquad \forall z\in\DD.
\end{equation}
\item[\rm (iii)] There exists 
a Hilbert space $\EE$, a conjugation $J'$ on $\EE$, and
a pure contractive analytic function $\Theta:\DD\to\LL(\EE)$, whose values are $J'$-symmetric operators, such that  $\Theta_T$ coincides with $\Theta$.
\end{enumerate}
\end{Thm}

\begin{proof}
(i)$\Rightarrow$(ii)\quad If $T$ is complex symmetric, there exists a conjugation $C$ on $\HH$ such that $T=CT^*C$. Since $C$ is involutive, we get $CT^*=TC$, $CT=T^*C$, and 
$C(I-T^*T)=(I-TT^*)C$. Thus $CD_T^2=D_{T^*}^2C$, and therefore $CD_T^{2n}=D_{T^*}^{2n}C$, $n\geq 0$. If $(p_n)_{n\geq 1}$ is a sequence of polynomials tending uniformly to $\sqrt x$ on $[0,1]$,  then 
$Cp_n(D_T^2)=p_n(D_{T^*}^2)C$, whence
$CD_T=D_{T^*}C$. In particular, $C\mathcal D_T\subset\mathcal D_{T^*}$; since $T^*$ is 
also $C$-symmetric, we actually have equality. 
Moreover, $CT^n={T^*}^nC$ for all $n\ge 1$ implies 
$C(I-\overline z T)^{-1}=(I-zT^*)^{-1}C.$

Define now $J:=C|\mathcal D_T$. Then  $J$ is an anti-linear map from $\mathcal D_T$ onto $\mathcal D_{T^*}$ which is isometric, and  the equalities above imply that 
$J\Theta_T(z)^*J=\Theta_T(z)$ for all $z\in\DD$.

(ii)$\Rightarrow$(i)\quad Assume first that $T$ is completely non-unitary. We will prove that the model operator $\mathbf T_{\Theta_T}\in\LL(\mathfrak K_{\Theta_T})$, as defined by
\eqref{eq:espace-modele1} and \eqref{eq:operateur-modele}, is complex symmetric. For simplicity,
we will write in the sequel of the proof $\mathbf T$ and $\mathfrak K$ instead of $\mathbf T_{\Theta_T}$ and $\mathfrak K_{\Theta_T}$.

Let us introduce some supplementary notations. Define
\[\mathfrak H:=L^2(\DTT)\oplus \overline{(I-\Theta_T^*\Theta_T)^{1/2} L^2(\DT)}\]
and $\pi:L^2(\DT)\to \mathfrak H$, $\pi_*:L^2(\DTT)\to \mathfrak H$ by
\[
\pi(f)=\Theta_T f\oplus (I-\Theta_T^*\Theta_T)^{1/2} f,\qquad\pi_*(g)=g\oplus 0,
\]
for $f\in L^2(\DT)$ and $g\in L^2(\DTT)$. Then $\pi$ and $\pi_*$ are isometries,
$\mathfrak H$ is spanned by $\pi L^2(\mathcal D_T)$ and $\pi_* L^2(\DTT)$, $\pi_*^*\pi=\Theta_T$,
and $\mathfrak K=\mathfrak H\ominus (\pi H^2(E)\oplus \pi_* H^2_-(E_*))$. If $P$ denotes the orthogonal projection (in $\mathfrak H$) onto $\mathfrak K$; then
$P=I_{\mathfrak H}-\pi P_+\pi^*-\pi_* P_-\pi_*^*$.

Let $Z\in\LL(\mathfrak H)$ be the unitary operator which acts as multiplication
by $z$ on both coordinates. Then $\pi(zf)=Z\pi f$, $\pi_*(zg)=Z\pi_*g$, and,
according to~\eqref{eq:operateur-modele}, $\mathbf T=PZ|\mathfrak K$.

If $\widetilde J:L^2(\DT)\to L^2(\DTT)$ is defined by 
$(\widetilde J f)(z)=\overline z J(f(z))$,  
then $\widetilde J$ is  anti-linear, isometric and onto; moreover
\begin{equation}\label{eq:prop-J}
\widetilde J P_+=P_-\widetilde J,\quad\widetilde J H^2(\DT)=H^2_-(\DTT),
\end{equation}
and $\widetilde J^{-1}g(z)=\bar z J^{-1}g(z)$ for $g\in L^2(\DTT)$.
 
We define the anti-linear map $C:\mathfrak H\to\mathfrak H$ by the formula
\[C\left(\pi f+\pi_* g \right):=\pi_*(\widetilde J f)+\pi({\widetilde J}^{-1} g),\qquad f\in L^2(\DT), g\in L^2(\DTT).\]

We  prove first that $C$ is a conjugation on $\mathfrak H$ and that $Z$ is $C$-symmetric.
Since $\pi, \pi_*,\widetilde J,{\widetilde J}^{-1}$ are (linear or antilinear) isometries, 
it follows that for all $f,h\in L^2(E)$ and all $g,k\in L^2(E_*)$, 
\begin{align*}
\langle C(\pi f+\pi_* g),C(\pi h+\pi_* k) \rangle&=\langle \pi({\widetilde J}^{-1} g),\pi({\widetilde J}^{-1} k)\rangle+\langle \pi_*(\widetilde J f),\pi_*(\widetilde J h) \rangle\\
&\qquad\qquad+\langle \pi({\widetilde J}^{-1} g),\pi_*(\widetilde J h)\rangle+\langle \pi_*(\widetilde J f),\pi({\widetilde J}^{-1} k)\rangle\\
&=\langle k,g\rangle+\langle h,f \rangle+\langle \Theta_T {\widetilde J}^{-1} g,\widetilde J h\rangle+\langle \widetilde J f,\Theta_T {\widetilde J}^{-1} k\rangle.
\end{align*}
But $J\Theta_T(z)^*J=\Theta_T(z)$  implies $\Theta_T{\widetilde J}^{-1}=\widetilde J\Theta_T^*$,
and therefore
\begin{align*}
\langle C(\pi f+\pi_* g),C(\pi h+\pi_* k) \rangle=&\langle k,g\rangle+\langle h,f \rangle+\langle \widetilde J \Theta_T^*g,\widetilde J h\rangle+\langle \widetilde J f,\widetilde J \Theta_T^*k\rangle\\
=&\langle k,g\rangle+\langle h,f \rangle+\langle h,\Theta_T^*g\rangle+\langle  \Theta_T^*k,f\rangle\\
=&\langle k,g\rangle+\langle h,f \rangle+\langle h,\pi^*\pi_* g\rangle+\langle  \pi^*\pi_* k,f\rangle\\
=&\langle \pi h+\pi_* k,\pi f+\pi_* g \rangle.
\end{align*}
Thus $C$ is a well-defined isometric anti-linear map. It follows immediately from the  definition that $C^2=I_{\mathfrak H}$ and thus $C$ is a conjugation on $\mathfrak H$.

If $f\in L^2(\DT)$, then
\[
\begin{split}
CZC(\pi( f))&=CZ\pi_*(\widetilde J f)
=C\pi_*(z\widetilde J f)
=C\pi_*(Jf)\\
&=\pi(\widetilde J^{-1} J f)=\pi (\bar z J^{-1} J f)
=\pi(\bar z f)= Z^*\pi(f).
\end{split}
\]
Similarly one proves that $CZC(\pi_*(g))=Z^*\pi_*(g)$ for $g\in L^2(\DTT)$, and therefore
$CZC=Z^*$; that is, $Z$ is $C$-symmetric.

By \eqref{eq:prop-J}, $C(\pi H^2(\DT))=\pi_*\widetilde J H^2(\DT)=\pi_* H^2_-(\DTT)$ and $C(\pi_* H^2_-(\DTT))=\pi H^2(\DT)$. Since $C$ is isometric, we have
\[
C\mathfrak K=C\mathfrak H\ominus C\left(\pi H^2(\DT)\oplus \pi_* H^2_-(\DTT)\right)=\mathfrak H\ominus\left(\pi H^2(\DT)\oplus \pi_* H^2_-(\DTT)\right)=\mathfrak K.
\]
Therefore the restriction $C'$ of $C$ to $\mathfrak K$ is a conjugation on $\mathfrak K$.
Since $C$ leaves $\mathfrak K$ and its orthogonal invariant, we have 
$C|\mathfrak K=PCP|\mathfrak K$ and $PC(I_{\mathfrak H}-P)=0$. Therefore
\[
\mathbf T=PZ|\mathfrak K=PCZ^*C|\mathfrak K=PCPZ^*PCP|\mathfrak K=C'\mathbf T^*C'.
\]
Thus $\mathbf T$ is $C'$-symmetric. Since $T$ is completely non-unitary, $T$ is unitarily equivalent to $\mathbf T$ and is therefore also complex symmetric.

Now, let $T\in\LL(\HH)$ be a general contraction satisfying condition (ii) in the statement of the theorem. If we decompose $T=T_0\oplus T_u$, with $T_0$ c.n.u. and $T_u$ unitary, then $T_0$ also satisfies (ii), and it is therefore complex symmetric by the above argument. Since $T_u$ is unitary, it is complex symmetric. Therefore $T$, being the direct sum of two complex symmetric operators,
is also complex symmetric.

(ii)$\Rightarrow$(iii)\quad If
\eqref{eq:caract-symmetric} is satisfied, and $C'$ is some conjugation on
$\DT$, then $U=JC':\DT\to\DTT$ is unitary and $C'=U^*J$. If
$\Theta:\DD\to\LL(\DT)$ is defined by $\Theta(z)=U^*\Theta_T(z)$, then
\[
\Theta(z)= U^* J\Theta_T(z)^*J= U^*J (U^*\Theta_T(z))^* U^*J=
C' \Theta(z)^* C'.
\]

(iii)$\Rightarrow$(ii)\quad If $U:\EE\to\DT$, $U_*:\EE\to \DTT$ are unitary
operators satisfying $\Theta_T(z)=U_*\Theta(z) U^*$ for all $z\in\DD$, then $J=U_*J'U^*$
satisfies all requirements in~(ii).
\end{proof}

\begin{Cor}\label{co:def1}
A contraction $T$ with $\partial_T=\partial_{T^*}=1$ is complex symmetric.
\end{Cor}

\begin{proof}
If $\partial_T=\partial_{T^*}=1$, then $\Theta_T$ is scalar-valued and we may identify $\DT$ and $\DTT$ with $\CC$. The natural conjugation $J$ on $\CC$ defined by $J(z)=\overline z$ satisfies then 
condition (iii) in Theorem~\ref{Thm:main}, whence $T$ is complex symmetric.  
\end{proof}

For the case $T\in C_{00}$, Corollary~\ref{co:def1} is proved in \cite{GP} and \cite{Garcia2},
where more of its consequences
are developed. Also in~\cite{GP} one can find the next result, for which we give a different
proof.

\begin{Cor}\label{co:def2}
Any operator on a 2-dimensional space is complex symmetric.
\end{Cor}

\begin{proof}
Since the complex symmetry is preserved by multiplication with non-zero scalars,
it is enough to assume $\|T\|=1$. But then either $T$ is unitary, or $\partial_T=\partial_{T^*}=1$,
in which case we may apply Corollary~\ref{co:def1}.
\end{proof}

It follows from Theorem~\ref{Thm:main} that if a contraction $T$ is complex symmetric,
then $\partial_T=\partial_{T^*}$. However, this is also a consequence of a more general
result from~\cite{GP}, namely that if a (not necessarily contractive) operator $T$ is complex symmetric, then $\dim\ker T=\dim\ker T^*$.

\section{$2\times2$ inner functions}

As shown in Corollary~\ref{co:def1}, contractions with defect indices 1 are always
complex symmetric.  As an application of
Theorem~\ref{Thm:main}, we will discuss in this section the case $\partial_T=\partial_{T^*}=2$.
We assume moreover that the characteristic function $\Theta_T$ is inner, which is 
equivalent to $T\in C_{00}$.

\begin{Def}
Let $\Theta:\DD\to\LL(\EE,\EE_*)$ be a contractive analytic function. We say that $\Theta$ is \emph{symmetrizable} if its matrix with respect to some fixed orthonormal bases (independent of $z$) in $\EE$ and $\EE_*$ is symmetric for all $z\in\DD$.
\end{Def}

According to Lemma~\ref{Prop:symmetrie} and Theorem~\ref{Thm:main}, (iii), a contraction 
is complex symmetric if and only if its characteristic function is symmetrizable. We
are interested in this section in $2\times2$-matrix valued characteristic
functions $\Theta(z)$. Note that Corollary~\ref{co:def2} implies that, for all $z\in\DD$, there exist $U_1(z),U_2(z)$ unitary such that $U_1(z)\Theta(z)U_2(z)$ is symmetric. But, in order to find
symmetrizable analytic functions, the matrices $U_1$ and $U_2$ should not depend on $z$.

We recall the following result in~\cite{Garcia1} which gives a parametrization of $2\times 2$ inner functions.

\begin{Prop}\label{Prop:inner}
Suppose $\phi$ be a non constant inner function in $H^\infty$, $a,b,c,d\in H^\infty$, and 
\[\Theta(z)=
\begin{pmatrix}
a(z)&-b(z) \\
c(z)&d(z)
\end{pmatrix},\]
Then $\Theta$ is a $2\times2$ inner function and $\hbox{det }\Theta=\phi$ if and only if
\begin{enumerate}
\item[$\mathrm{(i)}$] $a,b,c,d$ belong to $\mathcal H(z\phi)=H^2\ominus z\phi H^2$;
\item[$\mathrm{(ii)}$] $d=C(a)$ and $c=C(b)$;
\item[$\mathrm{(iii)}$] $|a|^2+|b|^2=1$ a.e. on $\TT$.
\end{enumerate}
Here $C$ denotes the natural conjugation on $\HH_{z\phi}$ defined by
\begin{equation}\label{eq:C}
C(f)=\overline{f}\phi,\qquad (f\in \HH_{z\phi}).
\end{equation}
\end{Prop}

The following result characterizes the symmetrizable $2\times2$-matrix valued inner functions.

\begin{Thm}\label{thm:caract-symetrisable}
A $2\times 2$  inner function
$\Theta(z)=\left(\begin{smallmatrix}
a(z) & -b(z) \\ \noalign{\smallskip}
C(b)(z)& C(a)(z)
\end{smallmatrix}\right)$
is symmetrizable if and only if there exist $(\gamma,\theta)\not=(0,0)$ such that $\gamma a +\theta b$ is
a fixed point of~$C$, where $C$ is defined by~\eqref{eq:C}, $\phi=\hbox{det }\Theta$.
\end{Thm}

%
%

\begin{proof} Suppose there exists $(\gamma,\theta)\not=(0,0)$ such that $C(\gamma a+\theta b)=
\gamma a+\theta b$; we may
 assume that  $|\gamma|^2+|\theta|^2=1$. Define the unitary matrix $U$ by 
$
U=\left(\begin{smallmatrix}
\overline{\theta} & -\gamma  \\
\overline{\gamma} & \theta
\end{smallmatrix}\right)$.
Then
\[
\begin{pmatrix}
-i&0\\
0&i
\end{pmatrix}
\Theta(z)U=\begin{pmatrix}
\overline{\theta}a(z)-\overline{\gamma}b(z)& i(\gamma a(z)+\theta b(z))\\
i(\gamma a(z)+\theta b(z))& -\gamma C(b)(z)+\theta C(a)(z)
\end{pmatrix}
\]
since $\overline\theta C(b)(z)+\overline\gamma C(a)(z)=C(\gamma a+\theta b)(z)=(\gamma a+\theta b)(z)$. Therefore $\Theta$ is symmetrizable.

Reciprocally, assume that $\Theta$ is symmetrizable. If
a nontrivial linear combination
of $a,b$ is 0, then we are done, since of course 0 is a fixed point of $C$.

Suppose then that the system $\{a,b\}$ is linearly independent.  By definition, there exists two unitary matrices $U_1$ and $U_2$ such that $U_1 \Theta(z) U_2$ is symmetric for all $z\in\DD$.
Write
\[
U_1=\begin{pmatrix}
\mu & -\overline\lambda \\
\lambda & \overline\mu
\end{pmatrix},\qquad U_2=\begin{pmatrix}
\theta & -\overline\gamma \\
\gamma & \overline\theta
\end{pmatrix},
\]
with $|\mu|^2+|\lambda|^2=1$ and $|\theta|^2+|\gamma|^2=1$. 
Straightforward computations show that
\[
U_1 \Theta(z) U_2=\begin{pmatrix}
\ast &X \\
Y & \ast
\end{pmatrix},
\]
with $X=-\mu(\overline\gamma a+\overline\theta b)-\overline\lambda C(-\gamma b+\theta a)$ and $Y=\lambda(\theta a-\gamma b)+\overline\mu C(\overline\theta b+\overline\gamma a)$. Then the symmetry of the matrix is equivalent to
$$-(\mu\overline\gamma+\lambda\theta)a-(\mu\theta-\lambda\gamma)b=C\left((\mu\overline\gamma+\lambda\theta)a+(\mu\overline\theta-\lambda\gamma)b\right).$$
If we put $u:=(\mu\overline\gamma+\lambda\theta)a+(\mu\overline\theta-\lambda\gamma)b$, then it follows that $C(u)=-u$, that is $C(iu)=iu$, and $iu$ is a fixed point of~$C$. To conclude 
the proof, we need to show that $(\mu\overline\gamma+\lambda\theta,\mu\overline\theta-\lambda\gamma )\not=(0,0)$.

Suppose then that
\begin{equation}\label{eq:111}
\mu\overline\gamma+\lambda\theta=\mu\overline\theta-\lambda\gamma=0.
\end{equation}
If we multiply $\mu\overline\gamma=-\lambda\theta$  by $\overline\theta$ and 
$\mu\overline\theta=\lambda\gamma$ by
$\overline\gamma$, 
and substract, we obtain
$\lambda(|\theta|^2+|\gamma|^2)=0.$ But $|\theta|^2+|\gamma|^2=1$, so $\lambda=0$, whence $|\mu|^2+|\lambda|^2=1$ yields $|\mu|=1$. Then~\eqref{eq:111} implies  $\gamma=\theta=0$: a contradiction. 
\end{proof}

\begin{remark}
Note that the fixed points of a conjugation $C$ can easily be described by using Lemma~\ref{Prop:symmetrie}~(i).
They form the real vector space of all elements which have real Fourier coefficients
with respect to a $C$-real orthonormal basis.
\end{remark}

\begin{remark}
A closely related question would be to describe all symmetric $2\times 2$-matrix valued
analytic contractive inner functions $\Theta(z)$. This can be done along the lines of the solution
of the Darlington sythesis problem in~\cite[Section 5]{Garcia1},
as follows. We fix first $\det \Theta$, which will be a nonconstant scalar inner function $\phi\in H^\infty$. Then we take a function $b\in H^2\ominus z\phi H^2$,
such that $Cb=b$ ($C$ the conjugation $f\mapsto \phi \bar f$ on $H^2\ominus z\phi H^2$).
If $b$ is inner, then $b^2=\phi$, and 
\[\Theta(z)=
\begin{pmatrix}
0 & ib(z)\\ ib(z) & 0
\end{pmatrix}.
\]
If $b$ is not inner, then we take $a\in H^2\ominus z\phi H^2$, such that $|a|^2+|b|^2=1$
(such $a$'s exist by~\cite[Proposition 5.2]{Garcia1}). Then 
\[
\Theta(z)=\begin{pmatrix}
a(z) & ib(z)\\ ib(z) & C(a)(z)
\end{pmatrix}.
\]

In~\cite[8.2]{GP} one discusses further the parametrization of all rational solutions of
a Darlington synthesis. Similarly, one could describe all rational symmetric $2\times 2$-matrix valued
analytic contractive inner functions $\Theta(z)$.

However, our interest is rather in complex symmetric contractions, and
the characteristic function is only a method of studying them.
A parametrization of all
complex symmetric contractions with defect indices~2 would require also to determine when
two characteristic functions as above coincide. This problem does not seem to have a neat solution.
\end{remark}

\section{An example}

Consider two nonconstant scalar inner functions $u,v\in H^\infty$, and let $\TTT_u, \TTT_v$ be the
corresponding model operators (the \emph{compressed shifts} in the terminology of~\cite{Garcia1}).
The corresponding model spaces are $\KK_u=H^2\ominus uH^2$ and $\KK_v=H^2\ominus vH^2$. As noted above, $\TTT_u$ and $\TTT_v$ are completely non unitary contractions with characteristic
functions $u$ and $v$ respectively. Their defect spaces are 1-dimensional, and it follows from Corollary~\ref{co:def1} that they are 
both complex symmetric.

We will discuss the contractions of the form
\begin{equation}\label{eq:ex1}
T=\begin{pmatrix} \TTT_u & X\\ 0 & \TTT_v
\end{pmatrix};
\end{equation}
thus $T\in\LL(\HH)$, with $\HH=\KK_u\oplus \KK_v$.
The next lemma gathers some facts about this operator. 

\begin{Lem}\label{le:ex1}  Suppose $T\in\LL(\KK_u\oplus\KK_v)$ is a contraction. Then:
\begin{enumerate}
\item[$\mathrm{(i)}$] $X=D_{\TTT_u^*}YD_{\TTT_v}$, with $Y:\Dcal_{\TTT_v}\to\Dcal_{\TTT_u^*}$ a contraction.
\item[$\mathrm{(ii)}$] $\partial_T=\partial_{T^*}=1$ if $\|Y\|=1$, and
$\partial_T=\partial_{T^*}=2$ otherwise.
\item[$\mathrm{(iii)}$] $T\in C_{00}$.
\end{enumerate}
\end{Lem}

Note that,  since $\dim\Dcal_{\TTT_v}=\dim\Dcal_{\TTT_u^*}=1$, $Y$ can actually be identified
with a complex number of modulus not larger than~1.

\begin{proof} The general form of the entries of a $2\times2$ contraction, as described,
for instance, in~\cite[Theorem 1.3]{AG} or~\cite[IV.3]{FF}, applied to the case when one
of the entries is null, yields immediately (i), as well as
an identification of $\Dcal_T$ with $\Dcal_{\TTT_v^*}\oplus\Dcal_Y$, and of
$\Dcal_{T^*}$ with $\Dcal_{\TTT_v}\oplus\Dcal_{Y^*}$, whence (ii) follows.

Finally, (iii) is an instance of a more general fact:
if $T=\left(\begin{smallmatrix}T_1 & X \\ 0 & T_2\end{smallmatrix}\right)$ is a contraction, then $T_i\in C_{0.}$ implies 
$T\in C_{0.}$. Indeed, take $\epsilon>0$, and a vector $x=x_1\oplus x_2$. Choose first 
$k$ such that $\|T_2^k x_2\|<\epsilon$. If $T^k(0\oplus x_2)= x'_1\oplus T_2^k x_2$,
take $k'$ such that $\|T_1^{k'} (x'_1+T_1^k x_1\|<\epsilon$. Then
\begin{align*}
\|T^{k+k'}x\|&= \|T^{k+k'}(x_1\oplus 0)+T^{k+k'}(0\oplus x_2)\|
=\|(T_1^{k+k'}x_1\oplus 0)+T^{k'} (x'_1\oplus T_2^k x_2)\|\\
&\le\|(T_1^{k'}(T_1^k x_1+x'_1)\oplus 0)\| +\| T^{k'}(0\oplus T_2^k x_2)\|\le 
\epsilon+\epsilon=2\epsilon.
\end{align*}
Since in our case $T_1=\TTT_u$ and $T_2=\TTT_v$ are both of class $C_{00}$, the result follows.
\end{proof}

The next theorem determines when is $T$ complex symmetric.

\begin{Thm}\label{th:example}
$T$ is complex symmetric precisely in the following cases:
\begin{enumerate} 
\item[$\mathrm{(i)}$] $Y=0$;
\item[$\mathrm{(ii)}$] $\|Y\|=1$;
\item[$\mathrm{(iii)}$] $0<\|Y\|<1$ and there exists $\lambda\in\DD$ and $\mu\in\TT$ such that $v=\mu b_\lambda(u)$, where $b_\lambda$ denotes the elementary Blaschke factor defined by
\[b_\lambda(z)=\frac {\lambda-z}{1-\overline\lambda z}.\]
\end{enumerate}
\end{Thm}

\begin{proof}
If $Y=0$, then $T=\TTT_u\oplus \TTT_v$, 
and is therefore complex symmetric as the direct sum of two complex symmetric
operators. 
If $\|Y\|=1$, then $T$ has defect indices~1 by Lemma~\ref{le:ex1}, and is therefore
symmetric by Corollary~\ref{co:def1}. (One can then see easily, using~\eqref{eq:fact} below, that $\Theta_T$ coincides
with the scalar function $uv$). We can thus suppose in the sequel that $0<\|Y\|<1$, and
$\partial_T=\partial_{T^*}=2$.
 
Since we intend to apply Theorem~\ref{Thm:main}, we have to determine
the characteristic function of~$T$. This can be calculated directly,
but in order to avoid some tedious computations, we prefer to
use the theory of invariant subspaces of contractions and factorizations
of the characteristic function, as developed in~\cite[Chapter VII]{SNF}. 

First, note that $T\in C_{00}$ implies $\Theta_T$  inner.
Since $\KK_u$ is an invariant subspace for $T$, it follows from Theorem~VII.1.1 and
Proposition VII.2.1. from~\cite{SNF} that one can factorize
\begin{equation}\label{eq:fact}
\Theta_T(z)=\Theta_2(z)\Theta_1(z)
\end{equation}
into two analytic inner functions, and that the characteristic functions of
$\TTT_u$ and $\TTT_v$, that is, $u$ and $v$, are equal to the pure parts of $\Theta_1$
and $\Theta_2$. Also, $\Theta_1$ and $\Theta_2$ being both inner, the dimensions
of their range spaces must both be equal to the dimension of the range of $\Theta_T$.

It follows then that $\Theta_1$
and $\Theta_2$ must be $2\times 2$ matrix valued inner functions, and their pure parts
are $u$ and $v$ respectively. They coincide therefore with
$\left(
\begin{smallmatrix} 1 & 0\\0 & u \end{smallmatrix}\right)$ and
$\left(\begin{smallmatrix} 1 & 0\\0 & v \end{smallmatrix}\right)$
respectively. According to~\eqref{eq:fact}, we have $2\times2$ unitary matrices $U_1,U_2,
V_1,V_2$ such that
\[
\Theta_T=U_1\begin{pmatrix} 1 & 0\\0 & u \end{pmatrix}U_2V_1
\begin{pmatrix} 1 & 0\\0 & v \end{pmatrix}
V_2,
\]
If we write 
\[
U_2V_1=\begin{pmatrix} \alpha & -\beta\\\bar\beta & \bar\alpha \end{pmatrix}
\]
with $\alpha$, $\beta$  complex numbers satisfying
$|\alpha|^2+|\beta|^2=1$,
it follows that the characteristic function $\Theta_T$ coincides with the inner function
\begin{equation}\label{eq:theta2}
\Theta(z)=\begin{pmatrix} 1 & 0\\0 & u \end{pmatrix}
\begin{pmatrix} \alpha & -\beta\\\bar\beta & \bar\alpha \end{pmatrix}
\begin{pmatrix} 1 & 0\\0 & v \end{pmatrix}=
\begin{pmatrix}\alpha & -\beta u(z) \\ \bar{\beta} v(z) & \bar{\alpha} u(z) v(z)\end{pmatrix}.
\end{equation}
Note that condition $0<\|Y\|<1$ implies both $\alpha$ and $\beta$ different from~0.

We apply now Theorem~\ref{thm:caract-symetrisable} in order to determine
when $\Theta$ as given by ~\eqref{eq:theta2} is symmetrizable. 
Since $\det \Theta=uv$, this happens
if and only if a linear combination of $\alpha$ and $\beta u$, not having both coefficients null,
belongs to the fixed points of the conjugation $C$ on $\KK_{zuv}$ given by $C(f)=uv\bar f$. 

If this is the case, and we write the combination as $g=s+tu$, $s,t\in\CC$ (and $s,t$ are not both null, which implies also $g\not=0$), then
\begin{equation}\label{eq:point-fixe}
C(g)=g\Leftrightarrow v(\overline s u+\overline t)=s+tu,
\end{equation}
and thus
\[
v=\frac{s+tu}{\overline s u+\overline t}.
\]
We must have $t\not=0$, since otherwise $uv$ is constant, which is not possible.
So we can write
$$v=\dfrac {t}{\overline t}\dfrac{\frac st+u}{1+\frac {\overline s}{\overline t}u}.$$
But now if $|s|=|t|$, then $v=\frac {t}{\overline s}$ which is impossible. If $|s|>|t|$, then we see that $v$ is at the same time analytic and coanalytic; whence $v$ is constant --- again a contradiction. So the only possibility is $|s|<|t|$. 
If we put $\lambda=-\frac {s}{t}$ and $\mu=-\frac {t}{\overline t}$ we get the desired conclusion that $v=\mu b_\lambda(u)$.

Conversely, suppose $v=\mu b_\lambda(u)$ with $|\lambda|<1$ and $|\mu|=1$. Write $\mu=-\frac{\zeta}{\overline\zeta}$, with $\zeta\not=0$. Then
\[
v=\frac{\zeta u-\lambda \zeta}{\overline\zeta-\overline\zeta\overline\lambda u},
\]
and if we define $s:=-\lambda\zeta$ and $t:=\zeta$, then 
$$v(\overline s u+\overline t)=v(-\overline\lambda\overline\zeta u+\overline\zeta)=\zeta u-\lambda\zeta=s+tu,$$
which implies by (\ref{eq:point-fixe}) that $C(g)=g$, with $g:=s+tu$. Since $t\not=0$, we may apply Theorem \ref{thm:caract-symetrisable} to conclude that $\Theta$ is symmetrizable. 

We have thus proved that in case $0<\|Y\|<1$, $\Theta_T$ is symmetrizable if and only if
$v=\mu b_\lambda(u)$ with $|\lambda|<1$ and $|\mu|=1$. Applying now Theorem~\ref{Thm:main}
ends the proof.
\end{proof}

Using Theorem~\ref{th:example}, it is easy to construct different examples
of complex symmetric and non complex symmetric operators with defect indices~2. 

It is not surprising that the condition obtained depends only on the norm
of $\|Y\|$ (or, rather, its modulus). Indeed, with a little effort one can
show that all operators $T$ corresponding to a fixed value of $\|Y\|$ are unitarily
equivalent.

\bibliography{cft_charsym}

\end{document}